\documentclass[12pt]{amsart}

\usepackage{amsmath,amssymb}

\newtheorem{thm}{Theorem}
\newtheorem{lem}[thm]{Lemma}
\newtheorem{prop}[thm]{Proposition}

\newtheorem{conj}[thm]{Conjecture}

\theoremstyle{plain}

\newcommand{\Z}{{\mathbf Z}}
\newcommand{\Q}{{\mathbf Q}}

\newcommand{\R}{{\mathbf R}}

\newcommand{\trace}{\operatorname{tr}}

\newcommand{\ideal}{I}
\newcommand{\T}{\TN}

\newcommand{\unit}{\alpha}

\newcommand{\weil}{U_N}
\newcommand{\sol}{\operatorname{Sol}}
\newcommand{\mnorm}{ {\mathcal N}} 
\newcommand{\bil}{B}
\newcommand{\diff}{L}

\newcommand{\TT}{\mathbf T^2}
\newcommand{\TN}{T_N}  
\newcommand{\OPN}{\operatorname{Op}_N}
\newcommand{\UN}{U_N}
\newcommand{\mat}{\operatorname{Mat}}

\newcommand{\OF}{\mathfrak{O}}
\newcommand{\Torusev}{C} 

\newcommand{\curve}{X} 

\newcommand{\torus}{C(N)} 
\newcommand{\torusN}{C(2N)}
\newcommand{\HN}{\mathcal H_N}

\newcommand{\re}{\operatorname{Re}}
\newcommand{\im}{\operatorname{Im}}

\newcommand{\sign}{\operatorname{sign}}

\newcommand{\tr}{\operatorname{tr}}

\newcommand{\fs}{f^\#}
\newcommand{\Ex}{\mathbf E}

\newcommand{\V}[1]{V_{#1}(\psi_j)}

\newcommand{\vnn}{\nu}
\newcommand{\vnm}{\mu}
\newcommand{\vnk}{\kappa}
\newcommand{\vnl}{\lambda}

\begin{document}

\title[Matrix elements for quantum cat maps]
{On  the distribution of matrix elements for the  
quantum cat map}

\date{January 10, 2004}

\begin{abstract}
For  many classically chaotic systems it is believed that the 
quantum wave functions become uniformly distributed,  
that is the matrix elements of smooth observables
tend to the phase space average of the observable. 
In this paper we study the fluctuations of the matrix elements 
for the desymmetrized quantum cat map.  
We present a conjecture for the distribution of the normalized
matrix elements, namely that their distribution is that of a certain
weighted sum of traces of independent matrices in $SU(2)$. 
This is in contrast to  generic chaotic systems where the distribution
is expected to be Gaussian. 
We compute the second and fourth moment of the normalized
matrix elements and obtain agreement with our conjecture. 
\end{abstract}

\thanks{
This work was supported in part by  the EC TMR network 
``Mathematical aspects of Quantum Chaos'' 
(HPRN-CT-2000-00103). 
P.K. was also supported in part by the NSF (DMS 0071503), the Royal Swedish
Academy of Sciences and the Swedish Research Council. 
Z.R. was also supported in part by the US-Israel Bi-National Science
Foundation.}

\author{P\"ar Kurlberg}
\address{
Department of Mathematics\\ 
Chalmers University of Technology\\
SE-412 96 Gothenburg  \\
Sweden}
\urladdr{www.math.chalmers.se/\~{ }kurlberg}
\email{kurlberg@math.chalmers.se}

\author{Ze\'ev Rudnick}
\address{Raymond and Beverly Sackler School of Mathematical Sciences,
Tel Aviv University, Tel Aviv 69978, Israel}
\email{rudnick@post.tau.ac.il}

\numberwithin{equation}{section}


\maketitle

\section{Introduction}

A fundamental feature of  quantum wave functions of classically
chaotic systems is that the matrix elements of smooth observables
tend to the phase space average of the observable, at least in the
sense  of convergence in the mean \cite{schnirelman-qe, CdV, Zelditch87} or
in the mean square \cite{Zelditch-c-star}. In many systems it is
believed that in fact {\em all} matrix elements converge to the
micro-canonical average, however this has only been demonstrated
for a couple of arithmetic systems: 
For ``quantum cat maps'' \cite{cat1}, 
and conditional on the
Generalized Riemann Hypothesis\footnote{An unconditional proof was
recently announced by Elon Lindenstrauss.}
also for the modular domain
\cite{Watson-thesis}, in both cases assuming that the systems are
desymmetrized by taking into account the action of ``Hecke operators''. 

As for the approach to the limit, it is expected that the  {\em
fluctuations} of the matrix elements about their limit are Gaussian
with variance given by classical correlations of the observable
\cite{feingold-peres-matrix-elements, EFKAMM}. 
In this note we study these fluctuations for the quantum cat map. 
Our finding is that for this system, the picture is very different.

We recall the basic setup \cite{BH, DE, DEGI, cat1} (see
section~\ref{background} for further background and any unexplained notation): 
The classical mechanical system is the iteration of a 
linear hyperbolic map $A\in SL(2,\Z)$  of the torus $\TT=\R^2/\Z^2$ 
(a ``cat map'').   
The quantum system is given by specifying an integer $N$, which 
plays the role of the inverse Planck constant. In what follows, $N$
will  be restricted to be a prime. 
The space of quantum states of the system is 
$\HN=L^2(\Z/N\Z)$.      
Let $f\in C^\infty(\TT)$ be a smooth, real valued observable and
$\OPN(f):\HN\to \HN$ its quantization.  
The quantization of the classical map $A$ is a unitary map $\UN(A)$ of
$\HN$. 

In \cite{cat1} we introduced {\em Hecke operators}, 
a group of commuting unitary maps of $\HN$, which 
commute with $\UN(A)$. The
space $\HN$ has an orthonormal basis consisting of  joint
eigenvectors $\{\psi_j\}_{j=1 }^N$ of $\UN(A)$, which we call {\em
Hecke  eigenfunctions}.    
The matrix elements $\langle\OPN(f)\psi_j,\psi_j\rangle$
converge\footnote{For arbitrary eigenfunctions, that is ones  which
are not Hecke eigenfunctions, this need not hold, see
\cite{scarred-cat}.} 
to the phase-space average $\int_{\TT} f(x) dx$ \cite{cat1}. 
Our goal is to understand their fluctuations around their limiting value. 

Our main result is to present a conjecture for the limiting  distribution 
of the normalized matrix elements 
$$
F_j^{(N)}:=\sqrt{N}\left(\langle \OPN(f)\psi_j,\psi_j\rangle -\int_{\TT}
f(x) dx \right)\;. 
$$
For this purpose, define a binary quadratic form associated to $A$
by
$$
Q(x,y) =cx^2+(d-a)xy-by^2,\qquad
A=\begin{pmatrix}a&b\\c&d\end{pmatrix}
$$
For an observable $f\in C^\infty(\TT)$ and an integer $\vnn$, set
$$
\fs(\vnn ):=\sum_{\substack {n=(n_1,n_2)\in \Z^2 \\ Q(n)=\vnn }} 
(-1)^{n_1 n_2} \^f(n)
$$
where $\^f(n)$ are the Fourier coefficients of $f$. 
\begin{conj}\label{conj1}
As $N\to\infty$ through primes, the limiting distribution of the
normalized matrix elements $F_j^{(N)}$ is that of the random variable
$$
X_f:= \sum_{\vnn \neq 0} \fs(\vnn) \tr(U_\vnn)
$$
where $U_{\vnn}$ are independently chosen random matrices in $SU(2)$
endowed with Haar probability measure.
\end{conj}

This conjecture predicts a radical departure from the Gaussian
fluctuations expected to hold for generic systems
\cite{feingold-peres-matrix-elements, EFKAMM}. Our 
first result confirms this conjecture for the variance of these
normalized matrix elements.
\begin{thm}
\label{thm:variance}
As $N\to \infty$ through primes, the variance of the normalized
matrix elements $F_j^{(N)}$ is given by
\begin{equation}\label{our variance} 
\frac 1N\sum_{j=1}^N |F_j^{(N)}|^2 \to \Ex(X_f^2) = \sum_{\vnn\neq 0}
|\fs(\vnn)|^2 \;.
\end{equation}
\end{thm}
For a comparison with the variance expected
for the case of {\em generic} systems, see Section~\ref{comparison}. 
A similar departure from this behaviour of the variance was
observed recently by Luo and Sarnak \cite{luo-sarnak-matrix-elements}
for the modular domain. For another analogy with that case, see
section~\ref{sec:diff}.

We also compute  the fourth moment of $F_j^{(N)}$ and
find agreement with Conjecture~\ref{conj1}:
\begin{thm}
\label{thm:fourth-moment}
The fourth moment of the normalized matrix elements is given by
$$
\frac 1N\sum_{j=1}^N |F_j^{(N)}|^4 \to \Ex(|X_f|^4) = 
2\sum_{\vnn\neq 0} |\fs(\vnn)|^4
$$
as $N\to \infty$ through primes.
\end{thm}


\noindent{\bf Acknowledgements:} We thank Peter Sarnak for discussions
on his work with Wenzhi Luo \cite{luo-sarnak-matrix-elements}.


\section{Background}\label{background}

The full details on the cat map and its quantization can be found in
\cite{cat1}. For the reader's convenience we briefly recall the setup:
The classical dynamics are given by a hyperbolic linear map
$A\in SL(2,\Z)$ so that
$x=(\begin{smallmatrix}p\\q\end{smallmatrix})\in \TT \mapsto Ax$ is a
symplectic map of the torus.  Given an observable $f\in
C^\infty(\TT)$, the classical evolution defined by $A$ is $f\mapsto
f\circ A$, where $(f\circ A)(x)=f(Ax)$. 

For doing quantum mechanics on the torus, one takes Planck's constant
to be $1/N$ and as the Hilbert space of
states  one takes $\HN:=L^2(\Z/N\Z)$,  
%
where the inner product 
is given by
\begin{equation*}
\langle \phi,\psi \rangle 
 = \frac1N \sum_{Q\bmod N} \phi(Q) \, \overline\psi(Q) .
\end{equation*}
 
The basic observables are given by the operators $\TN(n)$, $n\in
\Z^2$, acting on $\psi \in L^2(\Z/N\Z)$ via: 
\begin{equation}\label{action of T(n)}
\left( \TN(n_1,n_2)\psi \right) (Q) =
e^{\frac {i\pi n_1 n_2}N} e(\frac{n_2Q}N)\psi(Q+n_1).
\end{equation}
where 
$e(x) = e^{2\pi i x}$.  

For any smooth classical observable 
$f\in C^\infty(\TT)$ with  Fourier expansion 
$f(x) = \sum_{n\in \Z^2} \widehat f(n) e(nx)$, 
its quantization is given by 
$$
\OPN(f) := \sum_{n\in \Z^2} \widehat f(n)  \TN(n) \;.
$$

\subsection{Quantum Dynamics:}\label{dynamics} 
For $A$ which satisfies a certain parity condition, 
we can assign unitary operators $\UN(A)$,
acting on $L^2(\Z/N\Z)$, having the following important properties:
\begin{itemize}
\item 
``Exact Egorov'':  For all observables $f \in C^\infty(\TT)$
$$ 
 \UN(A)^{-1}  \OPN(f) \UN(A)= \OPN(f\circ A).
$$
\item The quantization depends only on $A$ modulo $2N$: If 
$A  \equiv B \mod 2N$ then  $\UN(A)=\UN(B)$. 
\item The quantization is multiplicative: if $A,B$ are 
congruent to the identity matrix modulo $4$ (resp., $2$) if $N$ is
even (resp., odd), then \cite{cat1, Mezzadri}
\begin{equation*}
\UN(AB)=\UN(A)\UN(B)
\end{equation*}
\end{itemize}

\subsection{Hecke eigenfunctions}\label{sec:Hecke eigenfunctions} 
Let $\unit$, $\unit^{-1}$ 
be the eigenvalues of $A$. 
Since $A$ is  
hyperbolic, $\unit$ is a unit in the real quadratic field $K=\Q(\unit)$.
Let $\OF=\Z[\unit]$, which is an order of $K$.
Let $v=(v_1,v_2)\in\OF^2$ be a vector such that 
$vA = \unit v$. 
If $A=\begin{pmatrix}a&b\\c&d\end{pmatrix}$, 
we may take $v=(c,\unit-a)$.
Let  $I:=\Z[v_1,v_2]=\Z[c,\unit-a]  \subset \OF$. 
Then $I$  is an $\OF$-ideal, and
the matrix of $\unit$ acting on $I$ by 
multiplication in the basis $v_1,v_2$ is precisely $A$. 
The choice of basis of $I$ gives an identification 
$I \cong \Z^2$ and the action of $\OF$  on the ideal $I$  by
multiplication gives  
a ring homomorphism
$$
\iota : \OF \to \mat_2(\Z)
$$ 
with  the property that $\det(\iota(\beta)) =\mnorm(\beta)$, where  
$\mnorm : \Q(\unit) \to \Q$ is the norm map. 

Let $\Torusev(2N)$ be the elements of $\OF/2N\OF$ with norm congruent
to $1\mod 2N$, and which congruent to $1$ modulo $4 \OF$ (resp., $2\OF$) if
$N$ is even (resp.,odd).  
Reducing $\iota$ modulo $2N$ gives a  map
$$
\iota_{2N} : \Torusev(2N) \to SL_2(\Z/2N\Z).
$$
Since $\Torusev(2N)$ is commutative, the multiplicativity of
our quantization 
implies that 
$$
\{ \UN( \iota_{2N}(\beta) ) : \beta \in \Torusev \}
$$
forms  a family of commuting operators.
Analogously with modular forms, we call these
{\em Hecke operators},
and functions $\psi \in \HN$ that are simultaneous eigenfunctions of
all the Hecke operators are denoted {\em Hecke eigenfunctions}. Note
that a Hecke eigenfunction is an eigenfunction of
$\UN(\iota_{2N}(\unit))=\UN(A)$.

The matrix elements are invariant under the Hecke operators:
$$
\langle\OPN(f)\psi_j,\psi_j \rangle  = 
\langle\OPN(f\circ B)\psi_j,\psi_j \rangle,\qquad B\in \torusN
$$
This follows from $\psi_j$ being eigenfunctions of the Hecke operators
$\torusN$. 
In particular, taking $f(x) = e(nx)$ we see that 
\begin{equation}\label{Hecke invariance}
\langle\TN(n)\psi_j,\psi_j \rangle = \langle\TN(nB)\psi_j,\psi_j
\rangle \;. 
\end{equation}

\subsection{The quadratic form associated to $A$: } 
We define a binary quadratic form associated to
$A=\begin{pmatrix}a&b\\c&d\end{pmatrix}$ by  
$$
Q(x,y) =cx^2+(d-a)xy-by^2
$$

This, up to sign, is the quadratic form 
$\mnorm(x c+y(\unit -a))/\mnorm(\ideal)$ induced by the norm form on the
ideal $\ideal = \Z[c,\unit-a]$ 
described in Section ~\ref{sec:Hecke eigenfunctions},    
where $\mnorm(I) = \#\OF/\ideal$. Indeed, since 
$\ideal=\Z[c,\unit-a]$ and $\OF=\Z[1,\unit]$ we have
$\mnorm(\ideal)=|c|$. A computation shows that the norm form is
then $\sign(c)Q(x,y)$. 

 By virtue of the  definition of $Q$ as a norm form, we see that $A$ and the
Hecke operators are isometries of $Q$, and since they have unit norm
they actually land in the special orthogonal group of $Q$.  That is we
find that under the above identifications, $\torusN $ is identified
with $\{B\in SO(Q,\Z/2N\Z): B\equiv I\mod 2 \}$.

\subsection{A rewriting of the matrix elements} 
We now show that when $\psi$ is a Hecke eigenfunction, the matrix 
elements $\langle\OPN(f)\psi,\psi \rangle$ have a modified Fourier
series expansion   which incorporates some extra invariance
properties.

\begin{lem}\label{lem:finite-support-equality}
If  $m,n\in \Z^2$ are such that
$Q(m)=Q(n)$, then for all sufficiently large primes $N$ we have 
$m\equiv nB \mod N$ for some $B\in SO(Q,\Z/N\Z)$. 
\end{lem}
\begin{proof}
We may clearly assume $Q(m)\neq 0$ because otherwise $m=n=0$ since $Q$
is anisotropic over the rationals. We take $N$ a sufficiently
large odd prime so that $Q$ is non-degenerate over the field $\Z/N\Z$. 
If $N> |Q(m)|$ then 
$Q(m)\neq 0\mod N$ and then the assertion reduces  to the 
fact that if $Q$ is a non-degenerate binary quadratic form 
over the finite field $\Z/N\Z$ ($N\neq 2$ prime)  then the 
special orthogonal group  $SO(Q,\Z/N\Z)$ acts transitively on the hyperbolas
$\{Q(n)=\vnn\}$, $\vnn\neq 0\mod N$. 
\end{proof}


%

\begin{lem}\label{lem:hyperbola} 
Fix $m,n\in \Z^2$ such that $Q(m)=Q(n)$. 
If $N$ is a sufficiently large odd prime and  $\psi$  a Hecke
eigenfunction, then 
$$
(-1)^{n_1 n_2}\langle \TN(n)\psi,\psi \rangle = 
(-1)^{m_1 m_2}\langle \TN(m)\psi,\psi \rangle 
$$
\end{lem}
\begin{proof}
For ease of notation, set $\epsilon(n):=(-1)^{n_1 n_2}$. 
By Lemma~\ref{lem:finite-support-equality} it suffices to show 
that if $m\equiv nB \mod N$ for
some $B\in SO(Q,\Z/N\Z)$ then 
$\epsilon(n)\langle \TN(n)\psi,\psi \rangle =  
\epsilon(m)\langle \TN(m)\psi,\psi \rangle$.  

By  the Chinese Remainder Theorem, 
$$
SO(Q,\Z/2N\Z)\simeq
SO(Q,\Z/N\Z)\times SO(Q,\Z/2\Z)
$$ 
(recall $N$ is odd) and so 
$$
\torusN\simeq\{B\in SO(Q\Z/2N\Z):B\equiv I\mod 2 \} \simeq 
SO(Q,\Z/N\Z)\times \{I\}
$$
Thus if $m\equiv nB \mod N$ for $B\in SO(Q,\Z/N\Z)$ then there is a
unique $\tilde B\in \torusN$ so that $m\equiv n\tilde B \mod N$.

We  note that $\epsilon(n)\TN(n)$ has period $N$, rather than
merely $2N$ for $\TN(n)$ as would follow from \eqref{action of T(n)}. 
Then since $m=n\tilde B \mod N$, 
$$
\epsilon(m)\TN(m) = \epsilon(n\tilde B) \TN(n\tilde B) =
\epsilon(n)\TN(n\tilde B)
$$
(recall that $\tilde B\in\torusN$ preserves parity: $n\tilde B\equiv
n\mod 2$, so $\epsilon(n\tilde B) = \epsilon(n)$). 
Thus for $\psi$ a Hecke eigenfunction, 
$$
\epsilon(m)\langle \TN(m)\psi,\psi\rangle = 
\epsilon(n)\langle \TN(n\tilde B)\psi,\psi\rangle = 
\epsilon(n)\langle \TN(n)\psi,\psi\rangle
$$
the last equality by \eqref{Hecke invariance}.  
\end{proof}

Define for $\vnn\in \Z$ 
$$
\fs(\vnn):= \sum_{n\in \Z^2: Q(n)=\vnn} (-1)^{n_1 n_2}\^f(n)
$$
and 
\begin{equation}
\label{eq:V-definition}
V_\vnn(\psi):= \sqrt{N} (-1)^{n_1n_2} \langle\TN(n)\psi,\psi \rangle
\end{equation}
where  $n\in \Z^2$ is a vector with $Q(n) = \vnn$ (if it exists)
and set $V_\vnn(\psi) = 0$ otherwise. By Lemma~\ref{lem:hyperbola} this
is well-defined, that is independent of the choice of $n$. 
Then we have 
\begin{prop}\label{prop:rewrite}
If $\psi$ is a Hecke eigenfunction, $f$ a trigonometric polynomial, 
and $N\geq N_0(f)$, then 
$$
\sqrt{N} \langle\OPN(f)\psi,\psi \rangle  = 
\sum_{\vnn\in \Z} \fs(\vnn) V_\vnn(\psi)
$$
\end{prop} 
To simplify the arguments, in what follows we will restrict ourself to
dealing with observables that are trigonometric polynomials.


\section{Ergodic averaging}
\label{sec:trace-trick}
We relate mixed moments of matrix coefficients to traces of certain
averages of the observables:  Let 
\begin{equation}
  \label{eq:definiton-of-D}
D(n)=
\frac{1}{|\torusN|}
\sum_{B \in \torusN}
\T(nB)
\end{equation}
The following shows that $D(n)$ is essentially diagonal when expressed in
the Hecke  eigenbasis.
\begin{lem}
\label{lem:almost-diagonal}
Let $\tilde{D}$ be the matrix obtained when expressing $D(n)$ in terms of
the 
Hecke eigenbasis $\{ \psi_i \}_{i=1}^N$.   If $N$ is inert in 
$K$, then $\tilde{D}$ is diagonal. 
If $N$ splits in $K$, then $\tilde{D}$ has the form
$$
\tilde{D}=
\begin{pmatrix}
D_{11} & D_{12} & 0 & 0 & \ldots & 0\\  
D_{21} & D_{22} & 0 & 0 & \ldots & 0\\  
0  & 0  & D_{33} & 0 & \ldots  & 0\\  
0 & 0  & 0  & D_{44} &  \ldots & 0\\  
\vdots & \vdots  & \vdots  & \vdots &  \ddots & \vdots\\  
0 & 0  & 0  & 0 &  \ldots & D_{NN}\\  
\end{pmatrix}
$$
where $\psi_1, \psi_2$ correspond to the quadratic character of
$\torusN$.  Moreover, in the split case, we have
$$|D_{ij}| \ll N^{-1/2}$$
for $1 \leq i,j \leq 2$.

\end{lem}
\begin{proof} 
  If $N$ is inert, then the Weil representation is multiplicity free
  when restricted to $\torusN$ (see Lemma 4 in \cite{lrh2}.)  If $N$ is
  split, then $\torusN$ is isomorphic to $(\Z/N\Z)^*$ 
and the trivial character occurs with multiplicity one, the quadratic
  character occurs with multiplicity two, and all other characters
  occur with multiplicity one (see \cite{catsup}, section 4.1). This
  explains the shape of $\tilde{D}$. 

As for the bound on in the split case, it suffices to take 
$f(x,y) = e( \frac{n_1x + n_2y}{N})$ for some 
$n_1,n_2 \in \Z$. 
%
We may give an explicit construction of the Hecke
eigenfunctions as follows (see \cite{catsup}, section 4 for more
details): there exists $M \in SL_2(\Z/2N\Z)$ such
that the eigenfunctions $\psi_1,\psi_2$ can be written as 
$$
\psi_1 = \sqrt{N} \cdot \weil(M) \delta_0,
\quad \psi_2 = \sqrt{\frac{N}{N-1}} \cdot  \weil(M) (1-\delta_0)
$$
where $\delta_0(x) = 1$ if $x \equiv 0 \mod N$, and $\delta_0(x) =
0$ otherwise.  
%
Setting  $\phi_1 = \sqrt{N} \delta_0$ and
$\phi_2 = \sqrt{\frac{N}{N-1}}(1-\delta_0)$,
exact Egorov gives
$$
D_{ij} = \langle \T((n_1,n_2)) \psi_i, \psi_j \rangle 
=
\langle \T((n_1',n_2')) \phi_i, \phi_j \rangle 
$$
where $(n_1',n_2') \equiv (n_1,n_2)M \mod N$.
Since we may assume $n$ not to be an eigenvector of $A$ modulo $N$, 
we  have $n_1'\not \equiv 0 \mod N$ and $n_2' \not \equiv 0 \mod N$.
Hence
$$
D_{11} =
\langle \T((n_1',n_2')) \phi_1, \phi_1 \rangle 
=
e(\frac{n_1'n_2'}{2N})\delta_0(0+n_1')
=
0
$$
since $n_1' \not \equiv 0 \mod N$.
The other estimates are analogous. 
%
%
\end{proof}

{\em Remark:}  In the split case, it is still true that $D_{ij} \ll
N^{-1/2}$ for {\em all} $i,j$, but this requires the Riemann
hypothesis for curves, whereas the above is elementary.

\begin{lem}
\label{lem:product-as-trace}
Let $\{\psi_i\}_{i=1}^N$ be a Hecke basis of $\HN$, and let $k,l,m,n
\in \Z^2$.  Then
$$
\sum_{i=1}^N
\langle \T(m) \psi_i, \psi_i \rangle  
\overline{\langle \T(n) \psi_i, \psi_i \rangle  }
=
\trace \big( D(m) D^*(n) \big)
+O(N^{-1})
$$
Moreover,
$$
\sum_{i=1}^N
\langle \T(k) \psi_i, \psi_i \rangle  
\overline{\langle \T(l) \psi_i, \psi_i \rangle  }
\langle \T(m) \psi_i, \psi_i \rangle  
\overline{\langle \T(n) \psi_i, \psi_i \rangle  }
$$
$$
=
\trace \big( D(k) D^*(l) D(m) D^*(n) \big)
+O(N^{-2})
$$
\end{lem}

\begin{proof}
By definition
$$
\sum_{i=1}^N
\langle \T(m) \psi_i, \psi_i \rangle  
\overline{\langle \T(n) \psi_i, \psi_i \rangle  }
= 
\sum_{i=1}^N
D(m)_{ii}
\overline{D(n)_{ii}}
$$
On the other hand, by lemma~\ref{lem:almost-diagonal},
$$
\trace \big( D(m) D(n)^* \big)
=
D_{12}(m) \overline{D_{21}(n)} +
D_{21}(m) \overline{D_{12}(n)} +
\sum_{i=1}^N
D_{ii}(m)
\overline{D_{ii}(n)}
$$
where $D_{12}(m), D_{21}(m), D_{12}(n)$ and $D_{21}(n)$ are all
$O(N^{-1/2})$.  Thus 
$$
\sum_{i=1}^N
\langle \T(m) \psi_i, \psi_i \rangle  
\overline{\langle \T(n) \psi_i, \psi_i \rangle  }
=
\trace \big( D(m) D(n)^* \big)
+O(N^{-1})
$$
The proof of the second assertion is similar.  
\end{proof}

\section{Proof of Theorem~\ref{thm:variance}}

In order to prove Theorem~\ref{thm:variance}
it suffices, by Proposition~\ref{prop:rewrite},
to show that as $N\to\infty$, 
$$
\frac{1}{N}
\sum_{j=1}^N
\V{\vnn}
\overline{\V{\vnm}}
\to 
\Ex \big(\trace U_\vnn \trace U_\vnm \big) 
=\begin{cases} 1& \text{ if } \vnm = \vnn, \\
0& \text{ if } \vnm\neq \vnn, \end{cases}
$$
where $U_\vnm, U_{\vnn} \in SU_2$ are  random
matrices in $SU_2$, independent if $\vnn \neq \vnm$.


\begin{prop}
\label{prop:V-variance}
Let $\{\psi_i\}_{i=1}^N$ be a Hecke basis of $\HN$.  
If $N \geq N_0(\vnm,\vnn)$ is prime and 
$\vnm,\vnn \not \equiv 0 \mod N$,
then 
$$
\frac{1}{N}
\sum_{j=1}^N
\V{\vnn}
\overline{\V{\vnm}}
=
\begin{cases}
1+ O(N^{-1}) & \text{if $\vnm = \vnn$,}\\
O(N^{-1}) & \text{otherwise.}
\end{cases}
$$  
\end{prop}
\begin{proof} 
Choose  $m,n \in \Z^2$  such that $Q(m)=\vnm$ and $Q(n)=\vnn$.
By (\ref{eq:V-definition}) and Lemma~\ref{lem:product-as-trace} we
find that 
\begin{equation*}
\begin{split}
\frac{1}{N}
\sum_{j=1}^N
\V{\vnn}
\overline{\V{\vnm}}
& = 
(-1)^{m_1m_2+n_1n_2} 
\sum_{j=1}^N
\langle \T(n) \psi_j, \psi_j \rangle  
\overline{\langle \T(m) \psi_j, \psi_j \rangle  }\\
& =
(-1)^{m_1m_2+n_1n_2} 
\trace \big( D(n) D(m)^* \big)
+O(N^{-1})
\end{split}
\end{equation*}
By definition of $D(n)$ we have 
$$
D(n) D(m)^* =
\frac{1}{|\torusN|^2}
\sum_{B_1,B_2 \in \torusN}
\T(nB_1)\T(mB_2)^* \;.
$$
We now take the trace of both sides and apply the following easily
checked identity (see \eqref{action of T(n)}), valid for 
odd $N$ and $B_1,B_2 \in \torusN$: 
$$
\trace(\T(nB_1) \T(mB_2)^*  ) =
\begin{cases}
(-1)^{m_1m_2+n_1n_2}N & \text{if $nB_1 \equiv mB_2 \mod N$,} \\
0 & \text{otherwise.}
\end{cases} 
$$
We get
\begin{multline}\label{ident-2}
\frac{1}{N}
\sum_{j=1}^N
\V{\vnn}
\overline{\V{\vnm}}
= \\
=
\frac{(-1)^{m_1 m_2+n_1n_2} }{|\torusN|^2}
\sum_{\substack{B_1,B_2 \in \torusN \\ nB_1 \equiv mB_2 \mod N}}
(-1)^{m_1 m_2+n_1n_2}N
+O(N^{-1}) \\
=
\frac{N }{|\torusN|} \cdot
|\{B \in \torusN : n \equiv mB \mod N\}|
+O(N^{-1})
\end{multline}
which, since $|\torusN| = N \pm 1$, equals
$1+ O(N^{-1})$ if there exists $B\in \torusN$ such that $ n
  \equiv mB \mod N$, and $O(N^{-1})$ otherwise.
Finally, for $N$ sufficiently large (i.e., $N \geq N_0(\vnm,\vnn)$),
Lemma~\ref{lem:finite-support-equality}
gives that $ n \equiv mB \mod N$ for some $B\in \torusN$ is equivalent
to $\vnm=\vnn$. 
\end{proof}


\section{Proof of theorem~\ref{thm:fourth-moment}}
\label{sec:proof-theorem-2}

\subsection{Reduction}\label{subsec:reduct}  
In order to prove Theorem~\ref{thm:fourth-moment} it suffices to show that 
\begin{equation}
\label{e:good-name}
\frac{1}{N}
\sum_{j=1}^N
\V{\vnk}
\overline{\V{\vnl}}
\V{\vnm}
\overline{\V{\vnn}}
\to 
\Ex \big(
\trace U_\vnk \trace U_\vnl 
\trace U_\vnm \trace U_\vnn \big) 
\end{equation}
where $U_\kappa, U_\lambda, U_\mu$ and $U_\nu$ are 
independent  random matrices in
$SU_2$.

Let $S \subset \Z^4$ be the set of four-tuples
$(\vnk,\vnl,\vnm,\vnn)$ such that $\vnk=\vnl,
\vnm=\vnn$, or $\vnk=\vnm, \vnl=\vnn$, or $\vnk=\vnn,
\vnl=\vnm$, but {\bf not} $\vnk= \vnl=\vnm=\vnn$.

\begin{prop}
\label{prop:V-fourth-moment}
Let $\{\psi_i\}_{i=1}^N$ be a Hecke basis of $\HN$ and let
$\vnk,\vnl,\vnm,\vnn \in \Z$.  If $N$ is a sufficiently large prime, 
then 
\begin{multline*}
\frac{1}{N}
\sum_{j=1}^N
\V{\vnk}
\overline{\V{\vnl}}
\V{\vnm}
\overline{\V{\vnn}}
 =
\begin{cases}
2+ O(N^{-1}) & \text{if $\vnk=\vnl =\vnm = \vnn $,}\\
1+ O(N^{-1}) & \text{if $(\vnk,\vnl,\vnm,\vnn) \in S$,}\\
 O(N^{-1/2}) & \text{otherwise.}\\
\end{cases}
\end{multline*}
\end{prop}

Given Proposition~\ref{prop:V-fourth-moment} it is straightforward to
deduce (\ref{e:good-name}), we need only to note that
$\Ex \big( (\trace U)^4 \big) =2$, 
$\Ex \big( (\trace U)^2 \big) =1$, and 
$\Ex \big( \trace U \big)=0$.

The proof of Proposition~\ref{prop:V-fourth-moment} will occupy
the remainder of this section. For the reader's convenience,  
here is a brief outline:
\begin{enumerate}
\item Express the left hand side of (\ref{e:good-name}) an exponential sum.
\item Show that the exponential sum is quite small unless pairwise
  equality of $\kappa,\lambda,\mu,\nu$ occurs, in which case the
  exponential sum is given by the number of solutions (modulo $N$) of
  a certain equation.
\item Determine the number of solutions.
\end{enumerate}

\subsection{Ergodic averaging} 
\begin{lem}\label{lem:tedious} 
Choose $k,l,m,n \in \Z^2$ such that $Q(k)
= \vnk, Q(l) = \vnl, Q(m) = \vnm$, and $Q(n) = \vnn$. Then 
\begin{multline}  \label{eq:10}
\frac{1}{N}
\sum_{j=1}^N
\V{\vnk}
\overline{\V{\vnl}}
\V{\vnm}
\overline{\V{\vnn}}
= 
\frac{N^2}{|\torusN|^4} 
\cdot \\ \cdot
\sum_{\substack{B_1, B_2, B_3, B_4 \in \torus \\
kB_1-lB_2+mB_3-nB_4 \equiv 0 \mod N}}
e \left( \frac{t(\omega(kB_1,-l B_2)+\omega(mB_3,-n B_4))  }{N} \right) 
\end{multline}
\end{lem} 
The proof of Lemma~\ref{lem:tedious} 
is an extension of the arguments proving the analogous
\eqref{ident-2} in the proof of Proposition~\ref{prop:V-variance}  
and is  left to the  reader.

\subsection{Exponential sums over curves}
\label{sec:nontrivial-exponential-sums}
In order to show that there is quite a bit of cancellation in 
(\ref{eq:10}) when pairwise
equality of norms do not hold, we will need some results on
exponential sums over curves.
Let $X$ be a projective curve of degree $d_1$ defined over the finite
field $\mathbb{F}_p$, embedded in $n$-dimensional projective space
$\mathbb{P}^n$ over 
$\mathbb{F}_p$.  Further, let $R(X_1, \ldots, X_{n+1})$ be a
homogeneous rational 
function in $\mathbb{P}^n$, defined over $\mathbb{F}_p$, and let $d_2$
be the degree of 
its numerator.  Define 
$$
S_m(R,X) = \sum_{x \in X(\mathbb{F}_{p^m})}' e\left( \frac{\sigma(R(x))}{p} 
\right)
$$
where $\sigma$ is the trace from $\mathbb{F}_{p^m}$ to $\mathbb{F}_p$, and  the
accent in the summation means that the poles of $R(x)$ are excluded.
%
\begin{thm}[Bombieri \cite{bombieri-exp-sum}, Theorem 6]
If $d_1d_2<p$ and $R$ is not constant on any component $\Gamma$ of
$X$ then 
$$
|S_m(R,X)|
\leq
(d_1^2+ 2d_1d_2 - 3d_1) p^{m/2} + d_1^2
$$
\end{thm}
In order to apply Bombieri's Theorem we need to show that the
components of a certain algebraic set are at most one
dimensional,
and in order to do this we show that the number of points defined over
$\mathbb{F}_N$ is $O(N)$. 
(Such a bound can not hold for all $N$ if there are
components of dimension two or higher.)
\begin{lem}
\label{l:two-circles}Let $a,b \in \mathbb{F}_N[\unit]$. 
If $a \neq 0$ and the equation 
$$
\gamma_1 = a \gamma_2 + b, \ \gamma_1,\gamma_2 \in \torus
$$
is satisfied for more than two values of $\gamma_2$, then $b =0 $ and
$\mnorm(a) = 1$.  
\end{lem}
\begin{proof}
Taking norms, we obtain $1 = \mnorm(a) + \mnorm(b) +
\trace(\overline{a} b\gamma_2)$ and 
hence $\trace(\overline{a} b\gamma_2)$ is constant.   If $\overline{a}
b \neq 0$, this  means that the coordinates $(x,y)$ of $\gamma_2$, when
regarding $\gamma_2$ as an element of $\mathbb{F}_N^2$, lies on some
line.  On the 
other hand, $\mnorm(\gamma_2)=1$ corresponds to $\gamma_2$ satisfying some
quadratic equation, hence the intersection can be  at most two points.
(In fact, we may identify $\torus$ with the solutions to
$x^2 -D y^2 = 1$ for $x,y \in \mathbb{F}_N$, and some fixed $D \in
\mathbb{F}_N$.) 
\end{proof}

\begin{lem}
\label{lem:why-curve}
Fix $k,l,m,n \in \Z^2$ and let $\curve$ be the set of solutions to 
$$
k-l B_2 + mB_3-n B_4 \equiv 0 \mod N, \ B_2,B_3,B_4 \in \torus
$$
If $Q(k), Q(l),Q(m),Q(n) \not \equiv 0 \mod N$, 
then $|\curve| \leq 3 (N+1)$ for $N$ sufficiently large.
\end{lem}

\begin{proof}
We use the identification of the action of $\torus$ on $\mathbb{F}_N^2$ with
the action of $\torus$ on $\mathbb{F}_N[\unit]$. 
The equation 
$$
k-l B_2 + mB_3-n B_4 \equiv 0 \mod N
$$
is then equivalent to
$$
\kappa - \lambda \beta_2 + \mu \beta_3 - \nu \beta_4 = 0
$$
where $\beta_i \in \torus$ and $\kappa,\lambda,\mu,\nu \in
\mathbb{F}_N[\unit]$.
We may rewrite this as 
$$
\kappa - \lambda \beta_2  =    \nu \beta_4 - \mu \beta_3=  
\beta_4 (\nu  - \mu \beta_3/\beta_4)
$$
and letting $\beta'=    \beta_3/\beta_4$, we obtain
$$
\kappa - \lambda \beta_2  =
\beta_4 (\nu  - \mu \beta') 
$$

If $\nu-\mu\beta'=0$ then $\kappa-\lambda\beta_2=0$, and since
$Q(l),Q(m) \not \equiv 0 \mod N$ implies that $\lambda,\mu$ are
nonzero\footnote{Recall that $Q$, up to a scalar multiple, is given by
  the norm.}, we find that $\beta_2$ and $\beta'$ are uniquely determined,
whereas $\beta_4$ can be chosen arbitrarily.  Thus there are at most
$|\torus|$ solutions for which $\nu-\mu\beta'=0$.

Let us now bound the number of solutions when $\nu-\mu\beta' \neq 0$:
after writing 
$$
\kappa - \lambda \beta_2  =
\beta_4 (\nu  - \mu \beta') 
$$
as 
$$
\frac{\kappa}{\nu  - \mu \beta'}
+ \frac{-\lambda}{\nu  - \mu \beta'} \beta_2  =
\beta_4,
$$
Lemma~\ref{l:two-circles} gives 
(note that $\kappa \neq 0$ since $Q(k) \not \equiv 0 \mod N$)
that there can be at most two
possible values of $\beta_2,\beta_4$ for each $\beta'$, and hence
there are at most $2|\torus|$ solutions for which $\nu-\mu\beta' \neq
0$. 
Thus, in total, $\curve$ can have at most $|\torus| + 2|\torus| \leq
3(N+1)$ solutions. 
\end{proof}

\subsection{Counting solutions}
\label{sec:counting-solutions}
We now determine the components of $X$ on which
$e \left( \frac{t(\omega(kB_1,-l B_2)+\omega(mB_3,-n B_4))  }{N}
\right)$ is constant.
\begin{lem}
\label{lem:number-of-diagonal-solutions}
Assume that $Q(k), Q(l),Q(m), Q(n) \not \equiv 0 \mod N$, and let
$\sol(k,l,m,n)$ be the number of solutions to the equations  
\begin{align}
\label{eq:5}
k B_1 - l B_2 + m B_3 -  n B_4  &\equiv 0 \mod N
\\
\label{eq:6}
\omega(kB_1,-l B_2)+\omega(mB_3,-n B_4) &\equiv -C \mod N
\end{align}
where $\ B_i \in \torus$.
%
%
If $C  \equiv 0 \mod N$ and 
$N$ is sufficiently large, then
\begin{equation}
  \label{eq:9}
\sol(k,l,m,n) =
\begin{cases}
2 |\torus|^2 & 
\text{if $Q(k) = Q(l) =Q(m) =  Q(n)$,} \\
|\torus|^2 + O(|\torus|) & 
\text{if $(Q(k),Q(l),Q(m),Q(n))  \in S$,}\\
O(|\torus|) &
\text{otherwise.}
\end{cases}
\end{equation}
On the other hand, if $C \not \equiv 0 \mod N$ then 
$$
\sol(k,l,m,n) = 
O(|\torus|).
$$
%
\end{lem}
\begin{proof}
  For simplicity\footnote{ The split case is similar except for
    possibility of zero divisors, but these do not occur when
    $k,l,m,n$ are fixed and  $N$ is
    large enough.}, we will assume that $N$
is inert. 
It will be convenient to use the language of algebraic number theory;
we identify $(\Z/N\Z)^2$ with the finite field
$\mathbb{F}_{N^2} = \mathbb{F}_N(\sqrt{D})$ by
letting $m = (x,y)$ correspond to $\mu=x+y\sqrt{D}$.  First we note
that if $n=  
(z,w)$ corresponds to $\nu$ then 
$$
\omega(m,n) = xw-zy = 
\im(  \overline{(x+y\sqrt{D})} ( z+w\sqrt{D} ))
$$
where  $\im(a+b\sqrt{D})= b$, and hence $\omega(m,n) =
\im(\overline{\mu}{\nu})$. 

Thus, with $(k,l,m,n)$ corresponding to $(\nu_1,\nu_2,\nu_3,\nu_4)$,
the values of $Q(k), Q(l), Q(m), Q(n)$ modulo $N$ are (up to a scalar
multiple) given by $\mnorm(\nu_1), \mnorm(\nu_2), \mnorm(\nu_3),
\mnorm(\nu_4)$. Putting $\mu_i = \nu_i \beta_i$ for $\beta_i \in
\torus$, we 
find that $\omega(kB_1,-lB_2) +\omega(mB_3,-nB_4) = -C$ can be written
as
$$
\im(\overline{\mu_1} \mu_2 + \overline{\mu_3} \mu_4) = C.
$$
Now, $k B_1 - l B_2 + m B_3 - n B_4 \equiv 0 \mod N$ is equivalent to
$\mu_1 - \mu_2 = \mu_4 -\mu_3$.  Taking norms, we obtain
$$
\mnorm(\mu_1) + \mnorm( \mu_2) - \trace(\overline{\mu_1} \mu_2)
=
\mnorm(\mu_4) + \mnorm( \mu_3) - \trace(\overline{\mu_4} \mu_3)
$$
and hence
$$
\trace(\overline{\mu_4} \mu_3)
= 
\trace(\overline{\mu_1} \mu_2)+ N_4+N_3-N_1-N_2
$$
if we let $N_i = \mnorm (\nu_i)$.
Since $\trace(\mu) = 2 \re(\mu) = 2 \re(\overline{\mu})$, we find that
$$
2\re(\overline{\mu_3} \mu_4)
= 
2\re(\mu_1 \overline{\mu_2})+ N_4+N_3-N_1-N_2
$$
On the other hand,  $\im(\overline{\mu_1} \mu_2 + \overline{\mu_3} \mu_4) = C$
implies that  
$$
\im(\overline{\mu_3} \mu_4)
= -\im(\overline{\mu_1} \mu_2) + C =
\im(\mu_1 \overline{\mu_2}) + C
$$
and thus
$$
\overline{\mu_3} \mu_4
=
\mu_1 \overline{\mu_2}
+K
$$
where $K = (N_4+N_3-N_1-N_2)/2 + C\sqrt{D}$.
Hence we can rewrite (\ref{eq:5}) and (\ref{eq:6}) as 
$$
\begin{cases}
\overline{\mu_3} \mu_4
=
\mu_1 \overline{\mu_2}
+K
\\
\mu_1 + \mu_3 = \mu_2 + \mu_4
\\
\mu_i = \nu_i \beta_i,  \ \beta_i \in \torus \text{ for $i=1,2,3,4$.}
\end{cases}
$$

\subsection*{Case 1 ($K \neq 0$)}
%
Since $\mu_i = \nu_i \beta_i$ with $\beta_i \in \torus$, we can rewrite 
$$
\overline{\mu_3} \mu_4
=
\mu_1 \overline{\mu_2}
+K
$$
as
$$
\overline{\nu_3} \nu_4  \beta_4/\beta_3
=
\nu_1 \overline{\nu_2} \beta_1/\beta_2
+K,
$$
and hence
$$
\beta_4/\beta_3
=
\frac{1}{\overline{\nu_3} \nu_4} 
( \nu_1 \overline{\nu_2} \beta_1/\beta_2+K).
$$
Applying lemma~\ref{l:two-circles} with $\gamma_1=\beta_4/\beta_3$ and
$\gamma_2=\beta_1/\beta_2$ gives that $\beta_1/\beta_2$, and hence
$\mu_1 \overline{\mu_2}$, must take one of two values, say $C_1$ or
$C_2$.
But $\mu_1 \overline{\mu_2} = C_1$ implies that 
$\mu_1 =  \mu_2 \frac{C_1}{N_2}$ and hence
$\mu_4 = \mu_3 \frac{C_1+K}{N_3} $.  We thus obtain
$$
\mu_2(1- \frac{C_1}{N_2}) = 
\mu_1-\mu_2=
\mu_4-\mu_3=
\mu_3(1- \frac{C_1+K}{N_3})
$$
Now, if $\mu_1 \neq \mu_2$ then both $1- \frac{C_1}{N_2}$ and $1-
\frac{C_1+K}{N_3}$ are nonzero.  Thus $\mu_2$ is determined by
$\mu_3$, which in turn gives that $\mu_1$ as well as $\mu_4$ are
determined by $\mu_3$.  Hence, there can be at most $\torus$ solutions
for which $\mu_1 \neq \mu_2$.  (The case $\mu_1 \overline{\mu_2} =
C_2$ is handled in the same way.)

On the other hand, for $\mu_1 = \mu_2$ we have the family of
solutions 
\begin{equation}
  \label{eq:diag-sol-two}
\mu_1 =\mu_2, \quad \mu_4=\mu_3
\end{equation}
(note that this implies that $C = 
\im(\overline{\mu_1} \mu_2+ \overline{\mu_3}\mu_4) = 0$.)

\subsection*{Case 2 ($K = 0$)}
%

Since $K=0$ and $\mu_1  = \mu_2 + \mu_4 - \mu_3$
we have
$$
\overline{\mu_3} \mu_4
=
\mu_1 \overline{\mu_2} + K
=
(\mu_2 + \mu_4 - \mu_3) \overline{\mu_2}
$$
and hence
$$
\mu_4 (\overline{\mu_3} -\overline{\mu_2})
=
(\mu_2  - \mu_3) \overline{\mu_2}
$$
If $\mu_2 - \mu_3 = 0$, we must have $\mu_1=\mu_4$, and we obtain
the family of solutions 
\begin{equation}
\label{eq:diag-sol-one}
\mu_2 = \mu_3,  \quad
\mu_1=\mu_4
\end{equation}

On the other hand, if $\mu_2  - \mu_3 \neq 0$, we can express 
$\mu_4$ in terms of $\mu_2$ and $\mu_3$:
$$
\mu_4 
=
\frac{\mu_2  - \mu_3}
{\overline{\mu_3} -\overline{\mu_2}}
\overline{\mu_2}
=
\frac{N_2  - \overline{\mu_2}\mu_3}
{N_3 - \overline{\mu_2}\mu_3}
\mu_3,
$$
which in turn gives that
\begin{multline}
  \label{eq:non-diag-sol}
\mu_1  = \mu_2 + \mu_4 - \mu_3
=
\mu_2 +
\frac{\mu_2  - \mu_3}
{\overline{\mu_3} -\overline{\mu_2}}
\overline{\mu_2}
- \mu_3
\\
=
\frac{\mu_2 - \mu_3}{\overline{\mu_3} -\overline{\mu_2}}
(\overline{\mu_3} -\overline{\mu_2})
+
\frac{\mu_2  - \mu_3}
{\overline{\mu_3} -\overline{\mu_2}}
\overline{\mu_2}
=
\frac{\mu_2  - \mu_3}
{\overline{\mu_3} -\overline{\mu_2}}
\overline{\mu_3}
=
\frac{\mu_2 \overline{\mu_3}  - N_3}
{\mu_2 \overline{\mu_3} - N_2}
\mu_2
\end{multline}

\subsection*{Summary}
If $K \neq 0$  there can be at most $2|\torus|$ ``spurious'' solutions
for which  $\mu_1 \neq \mu_2$; other than that, we must have
$$
\mu_1 = \mu_2,  \quad
\mu_3=\mu_4.
$$
On the other hand, if $K=0$, then either 
$$
\mu_2 = \mu_3,  \quad
\mu_1=\mu_4.
$$
or 
$$
\mu_4 
=
\frac{\mu_2  - \mu_3}
{\overline{\mu_3} -\overline{\mu_2}}
\overline{\mu_2}
=
\frac{N_2  - \overline{\mu_2}\mu_3}
{N_3 - \overline{\mu_2}\mu_3}
\mu_3,
\quad
\mu_1 =
\frac{\mu_2  - \mu_3}
{\overline{\mu_3} -\overline{\mu_2}}
\overline{\mu_3}
=
\frac{\mu_2 \overline{\mu_3}  - N_3}
{\mu_2 \overline{\mu_3} - N_2}
\mu_2
$$
We note that the first case can only happen if $N_1=N_2$ and
$N_3=N_4$, the second only if $N_2=N_3$ and $N_1=N_4$, and the third
only if $N_2=N_4$ and $N_1=N_3$.  Moreover, in all three cases,
$C=\im(K)= \im(\overline{\mu_1} \mu_2 + \overline{\mu_3}\mu_4 
) = 0$.  We also note that if $N_2=N_3$, then 
the third case simplifies to $\mu_1=\mu_2$ and
$\mu_3=\mu_4$.  We thus obtain the following:

If $C\neq 0$ then $K\neq 0$ and there can be at most $O(N)$
``spurious solutions''.  

If $C=0$ and $N_1=N_2=N_3=N_4$ then $K=0$ and the solutions are given
by the two families
$$
\mu_2 = \mu_3,  \quad
\mu_1=\mu_4
$$
and 
$$
\mu_4 
=
\frac{N_2  - \overline{\mu_2}\mu_3}
{N_3 - \overline{\mu_2}\mu_3}
\mu_3
= \mu_3,
\quad
\mu_1 
=
\frac{\mu_2 \overline{\mu_3}  - N_3}
{\mu_2 \overline{\mu_3} - N_2}
\mu_2
=
\mu_2
$$

If $C=0$ and $N_1=N_4 \neq N_2=N_3$ then $K=0$ and there is a family
of solutions given by 
$$
\mu_2 = \mu_3,  \quad
\mu_1=\mu_4.
$$

Similarly, if $C=0$ and $N_1=N_3 \neq N_2=N_4$ then $K=0$ and 
there is a family of solutions given by 
$$
\mu_4 
=
\frac{\mu_2  - \mu_3}
{\overline{\mu_3} -\overline{\mu_2}}
\overline{\mu_2},
\quad
\mu_1 =
\frac{\mu_2  - \mu_3}
{\overline{\mu_3} -\overline{\mu_2}}
\overline{\mu_3}
$$

If $C=0$ and $N_1=N_2 \neq N_3=N_4$ then $K \neq 0$, in which
case we have a family of solutions given by
$$
\mu_1= \mu_2, \quad \mu_3 = \mu_4
$$
as well as $O(N)$ ``spurious'' solutions.

Finally, if $C=0$ and pairwise equality of norms {\em do not} hold,
then we must have $K \neq 0$
(if $K=0$ then $\overline{\mu_3} \mu_4 = \mu_1 \overline{\mu_2}+K$ 
implies that $N_3N_4 = N_1N_2$, which together with $N_1+N_2 =
N_3+N_4$ gives that either $N_1=N_3, N_2=N_4$ or $N_1=N_4, N_2=N_3$)
and in this case there can be at most $O(N)$ ``spurious''
solutions.

Now Lemma~\ref{lem:finite-support-equality} gives 
that pairwise equality of norms modulo $N$
implies pairwise equality of $Q(k), Q(l),Q(m), Q(n)$.
\end{proof}


\subsection{Conclusion}
\label{sec:conclusion}
We may now evaluate the exponential sum in (\ref{eq:10})
\begin{prop} 
\label{prop:fourth-moment-expsum}
If $Q(k), Q(l), Q(m),  Q(n)  \not \equiv 0 \mod N$ then, for $N$
sufficiently large, we have
\begin{multline}
\label{eq:8}
\sum_{\substack{B_1, B_2, B_3, B_4 \in \torus \\
kB_1-lB_2+mB_3-nB_4 \equiv 0 \mod N}}
e \left( \frac{t(\omega(kB_1,-l B_2)+\omega(mB_3,-n B_4))  }{N} \right) 
\\
= \begin{cases}
2 |\torus|^2 + O(|\torus|) & \text{if $Q(k) = Q(l) =  Q(m) = Q(n)$, } \\
|\torus|^2 + O(|\torus|) & \text{if $(Q(k),Q(l),Q(m),Q(n)) \in S$,} \\
O(|\torus|^{3/2}) & \text{otherwise.} 
\end{cases}
\end{multline}
\end{prop}

\begin{proof}
Since both $\omega(kB_1,-l B_2)+\omega(mB_3,-n B_4)$ and $k B_1 - l
B_2 + m B_3 - n B_4$ are invariant under the substitution  
$
(B_1, B_2,B_3, B_4) \to (B' B_1, B' B_2 , B' B_3, B' B_4) 
$
for $B' \in \torus$,
we may rewrite the left hand side of (\ref{eq:8}) as $|\torus|$ times 
\begin{equation}
  \label{eq:expsum-answer}
\sum_{\substack{B_2, B_3, B_4 \in \torus \\
k-lB_2+mB_3-nB_4 \equiv 0 \mod N}}
e \left( \frac{t(\omega(k,-l B_2)+\omega(mB_3,-n B_4))  }{N} \right).
\end{equation}
Let $\curve$ be the set of solutions to 
$$
k-lB_2+mB_3-nB_4 \equiv 0 \mod N, \ B_2, B_3, B_4 \in \torus.
$$
By Lemma~\ref{lem:why-curve}, the dimension of any irreducible
component of $X$ is at most $1$.  The contribution
from the zero dimensional components of $\curve$ is at most
$O(|\torus|)$. 
As for the one dimensional components, 
Lemma~\ref{lem:number-of-diagonal-solutions} gives that $\omega(k,-l
B_2)+\omega(mB_3,-n B_4)$ cannot be constant on any  
component unless pairwise equality of norms holds.
%
Thus, 
if pairwise equality of norms does {\em not} hold, Bombieri's Theorem
gives that \eqref{eq:expsum-answer} is $O(N^{1/2}) =
O(|\torus|^{1/2})$. 

On the other hand, if
$\omega(kB_1,-l B_2)+\omega(mB_3,-n
B_4)$ equals some constant $C$ modulo $N$ on some one dimensional
component, then Lemma~\ref{lem:number-of-diagonal-solutions} 
gives the following:
$C \equiv 0 \mod N$, 
and (\ref{eq:expsum-answer}) equals $\sol(k,l,m,n)$, which in turn
equals $|\torus|^2$ or $2 |\torus|^2 $ depending on whether $Q(k) \equiv
Q(l) \equiv Q(m) \equiv Q(n) \mod N$ or not.
\end{proof}

Proposition~\ref{prop:V-fourth-moment} now follows from
 Lemma~\ref{lem:tedious} and 
Proposition~\ref{prop:fourth-moment-expsum} on recalling that 
$|\torus|= |\torusN|=N\pm1$.


\section{Discussion} 

\subsection{Comparison with generic systems} \label{comparison}
It is interesting to compare our result for the variance with 
the predicted answer for generic systems 
(see \cite{feingold-peres-matrix-elements, EFKAMM}), which is 
\begin{equation}\label{generic answer} 
\sum_{t=-\infty}^\infty \int_{\TT}  f_0(x) \overline{f_0(A^tx)} dx
\end{equation}
where  $f_0=f-\int_{\TT} f(y)dy $. Using the Fourier expansion
and collecting together frequencies $n$ lying in the same $A$-orbit
this equals 
$$
\sum_{t=-\infty}^\infty \sum_{0\neq n\in \Z^2} 
\^f(n)\overline{\^f(n A^t)} 
=
\sum_{m\in (\Z^2-0)/\langle A \rangle} \left| \sum_{n\in m\langle A
\rangle} \^f(n) \right|^2 
$$
where $\langle A \rangle$ denotes the group generated by $A$. We can
further rewrite this expression into a form closer to our 
formula \eqref{our variance} by noticing that the expression
$\epsilon(n):=(-1)^{n_1n_2}$ is an invariant of the $A$-orbit:
$\epsilon(n)=\epsilon(nA)$, because we assume that $A\equiv I \mod
2$. Thus we can write the generic variance \eqref{generic answer} as 
\begin{equation}\label{generic variance}
\sum_{m\in (\Z^2-0)/\langle A \rangle} \left| \sum_{n\in m\langle A
\rangle} (-1)^{n_1 n_2}\^f(n) \right|^2 \;.
\end{equation}
The comparison with with our answer 
$
\sum_{\vnn\neq 0} \left|\sum_{Q(n)=\vnn} (-1)^{n_1n_2} \^f(n) \right|^2  
$
in \eqref{our variance}, is now clear: Both expressions would coincide if each hyperbola
$\{n\in \Z^2 : Q(n)=\vnn\}$ consisted of a single $A$-orbit. It is
true that each 
hyperbola consists of a finite number of $A$-orbits for $\vnn\neq 0$, 
but that number varies with $\vnn$. 

\subsection{A differential operator}\label{sec:diff}
There is yet another analogy with the modular domain, pointed out to
us by Peter Sarnak: 
We define a differential operator $\diff$ on $C^\infty(\TT)$ by
$$
\diff = -\frac 1{4\pi^2} 
Q( \frac{\partial}{ \partial p}, \frac{\partial}{ \partial q})
$$
so that $\widehat{\diff f}(n) = Q(n)\^f (n)$. 

Given observables $f,g$, we define a bilinear form $\bil(f,g)$ by
$$
\bil(f,g)
= \sum_{\vnn\neq 0} \fs(\vnn) g^\#(\vnn)
$$
so that (cf. Conjecture~\ref{conj1}) 
$
\bil(f,g) = \Ex(X_f X_g)
$
and by Theorem~\ref{thm:variance},   
$B(f,f)$ is the variance of the normalized matrix elements.

It is easy to check that $\diff$ is self adjoint with respect to $\bil$, i.e., 
$\bil(\diff f, g) = \bil(f, \diff g)$. 
Note that $\diff$ is also self-adjoint with respect to the bilinear
form derived from the expected variance for generic systems
\eqref{generic answer}, \eqref{generic variance}.  This feature
was first observed for the modular domain, where the role of $\diff$
is played by the Casimir operator \cite{luo-sarnak-matrix-elements}
(c.f. Appendix 5 of Sarnak's survey
\cite{sarnak-spectra-of-hyperbolic-surfaces}).

\subsection{Connection with character sums} \label{sec:exp sums}
Conjecture~\ref{conj1}  is related to 
the value distributions of certain character sums, at least in 
the case of {\em split} primes, that is primes $N$ for which
the cat map $A$ is diagonalizable modulo $N$. 
%
%
Let $M \in SL_2(\Z/2N\Z)$ be such that 
$A=MDM^{-1}\mod 2N$.  
In \cite{catsup} we explained that in that case, all but
one 
of the normalized Hecke eigenfunctions are
given in terms of the Dirichlet characters $\chi$ modulo $N$ as 
$\psi_\chi:=\sqrt{\frac{N}{N-1}} \UN(M)\chi$.
We can then write the matrix elements 
$\langle\TN(n) \psi_\chi,\psi_\chi\rangle$ as characters sums: Setting 
$(m_1,m_2)=nM$, we have
$$
\langle\TN(n) \psi_\chi,\psi_\chi \rangle = 
e^{\pi i m_1 m_2/N} \frac 1{N-1} \sum_{Q\mod N} e(\frac {m_2 Q}N) 
\chi(Q+m_1)\overline{ \chi(Q)},
$$
%
and Conjecture~\ref{conj1} gives a prediction for the value
distribution of these 
sums as $\chi$ varies.



\bibliographystyle{abbrv} 

\end{document}